\documentclass[a4paper, 10pt]{amsart}

\usepackage{amsthm, amsmath, amssymb, amsfonts}
\usepackage{xcolor}
\usepackage{enumitem}
\usepackage{geometry}
\usepackage{mathrsfs}

\newgeometry{
	top=4cm,
	bottom=4cm,
	inner=3cm,
	outer=3cm,
}

\usepackage[colorlinks, linkcolor=cyan, citecolor=cyan, urlcolor=cyan, bookmarks=false, hypertexnames=true]{hyperref}

\theoremstyle{plain}
\newtheorem{theorem}{Theorem}[section]
\newtheorem{proposition}[theorem]{Proposition}

\newtheorem{corollary}[theorem]{Corollary}
\newtheorem{lemma}[theorem]{Lemma}

\newtheorem{innercustomgeneric}{\customgenericname}
\providecommand{\customgenericname}{}
\newcommand{\newcustomtheorem}[2]{%
	\newenvironment{#1}[1]
	{%
		\renewcommand\customgenericname{#2}%
		\renewcommand\theinnercustomgeneric{##1}%
		\innercustomgeneric
	}
	{\endinnercustomgeneric}
}

\newcustomtheorem{customthm}{Theorem}

\theoremstyle{definition}

\newtheorem{remark}[theorem]{Remark}

\numberwithin{equation}{section}

\DeclareMathOperator{\trace}{trace}
\DeclareMathOperator{\grad}{grad}

\DeclareMathOperator{\Div}{div}
\DeclareMathOperator{\id}{Id}

\DeclareMathOperator{\Riem}{Riem}

\allowdisplaybreaks

\title[]{Pinching rigidity results for compact biconservative hypersurfaces}

\author{\c Stefan Andronic, Aykut Kayhan and R\" uya Ye\u gin \c Sen}

\address{Faculty of Mathematics, Al. I. Cuza University of Iasi, Blvd. Carol I, no. 11, 700506 Iasi, Romania} \email{stefanandronic215@gmail.com, stefan.andronic@uaic.ro}

\address{Department of Mathematics, Near East University, Nicosia, Mersin-10, T\" urkiye} \email{aykut.kayhan@neu.edu.tr}

\address{Department of Mathematics, Faculty of Engineering and Natural Sciences, \. Istanbul Medeniyet University, 34700 \. Istanbul, \" Usk\" udar, T\" urkiye} \email{ruya.yegin@medeniyet.edu.tr}


\thanks{The second and third authors were supported by the Scientific and Technological Research Council of T\" urkiye (T\" UB\. ITAK), under the 3501-Career Development Program (CAREER) (project no: 125F553). The opinions and views expressed herein are those of the authors and do not reflect those of T\" UB\. ITAK}

\begin{document}
	\begin{abstract}
		It is known that the Nomizu--Smyth rigidity result for compact hypersurfaces with constant mean curvature and non-negative sectional curvature in space forms can be generalized by replacing the hypothesis of constant mean curvature with biconservativity, when the dimension of the hypersurface is at most $10$. In this paper, we extend this result to compact biconservative hypersurfaces of dimension greater than $10$ by imposing a pinching condition on the squared norm of the shape operator in terms of the mean curvature function. To obtain this result, we introduce a new family of symmetric $(1,1)$ tensor fields whose vanishing divergence characterizes the biconservativity of the hypersurface, and we use them in a tensorial formula which holds for any symmetric $(1,1)$ tensor field. When the mean curvature function does not vanish, we obtain a larger upper bound for the squared norm of the shape operator. We also obtain local versions of these results, in which the compactness of the hypersurface is replaced by the condition that a certain combination of the squared norm of the shape operator and the square of the mean curvature function is constant.
	\end{abstract}

	\maketitle
	
	\section{Introduction}
	
	Biharmonic maps $\varphi : \left ( M^m, g_M \right ) \to \left ( N^n, g_N \right )$ between Riemannian manifolds are a natural generalization of the well-known harmonic maps which are characterized by the vanishing of the tension field $\tau(\varphi) = \trace_{g_M} \nabla ^\varphi d \varphi$. We recall that, when $M$ is compact, a biharmonic map is a critical point of the bienergy functional 
	$$
	E_2 (\varphi) = \frac 1 2 \int _M |\tau (\varphi)|^2 \ v_g.
	$$ 
	Using the first variational formula, G. Y. Jiang proved that a biharmonic map is characterized by the vanishing of its bitension field
	\begin{equation*}\label{eq:Tau2Immersions}
		\tau_2(\varphi) = - \Delta ^\varphi \tau (\varphi) - \trace R^N \left ( d\varphi (\cdot), \tau (\varphi) \right ) d\varphi (\cdot),
	\end{equation*}
	where $\Delta^\varphi$ is the rough Laplacian acting on the sections of the pull-back bundle $\varphi^{-1} (TN)$ and $R^N$ is the curvature tensor field of the target manifold (see \cite{Jiang1986-1} and \cite{Jiang1986-2}). Since every harmonic map is biharmonic, we are interested in biharmonic maps that are not harmonic, which are called proper biharmonic. In the particular case of isometric immersions, that is, of submanifolds, the biharmonic equation $\tau_2 (\varphi) = 0$ splits into a tangent part and a normal part.
	
	Until now, all known examples of biharmonic submanifolds in space forms of non-positive sectional curvature are minimal, while in positively curved space forms only two classes of proper biharmonic submanifolds are known (see \cite{CaddeoMontaldoOniciuc2002} and \cite{Jiang1986-2}). Because of the scarcity of examples of proper biharmonic submanifolds in space forms, it was conjectured in \cite{Chen1991} that every biharmonic submanifold in the Euclidean space is minimal. In \cite{BalmusMontaldoOniciuc2008}, it was conjectured that the proper biharmonic submanifolds in the Euclidean sphere must have constant mean curvature (CMC).
	
	To overcome this scarcity, biconservative submanifolds were introduced in \cite{CaddeoMontaldoOniciucPiu2014}. They are characterized only by the vanishing of the tangent part of the bitension field. Besides having many more examples than biharmonic submanifolds (see, for example, \cite{AndronicMontaldoOniciucSanna2026}, \cite{FuHongYangZhan2025}, \cite{MontaldoOniciucPampano2023}, \cite{NistorOniciucTurgaySen2023}), biconservative submanifolds also have many interesting geometric properties, especially in the case of surfaces (see, for example, \cite{Nistor2017}). In the case of hypersurfaces, the biconservativity is characterized by the fact that the gradient of the mean curvature function $f$ and its image through the shape operator are parallel, with the proportionality factor $- m f / 2$. This implies that any CMC hypersurface is biconservative.
	
	It turns out that biconservativity is a natural notion in differential geometry, since it makes it possible to generalize some well-known results. One of them is the famous Nomizu--Smyth rigidity result which states that any compact CMC hypersurface $M^m$ with non-negative sectional curvature in a space form $N^{m+1} (c)$ has parallel shape operator (see \cite{NomizuSmyth1969} and Theorem \ref{th:NomizuSmythTheorem}). This result was generalized in \cite{FetcuOniciuc2022} by replacing the hypothesis of CMC with biconservativity when $m \leq 10$, and the same conclusion was obtained (see Theorem \ref{th:FetcuOniciuc2022}). Another example is the Cheng--Yau rigidity result (see \cite{ChengYau1977}) which was generalized to biconservative hypersurfaces in \cite{FetcuLoubeauOniciuc2021}. Shorter and simpler proofs of these generalizations were given in \cite{AndronicKayhan2025}, together with other rigidity results for compact biconservative hypersurfaces.
	
	In \cite{FetcuOniciuc2022}, the generalization of the Nomizu--Smyth rigidity result was obtained under a condition on the dimension of the hypersurface, namely $m \leq 10$. In this paper, we extend this result to compact biconservative hypersurfaces of dimensions greater than $10$ by imposing a pinching condition on the squared norm of the shape operator in terms of the mean curvature function.
	
	The main difficulty with the techniques used in \cite{AndronicKayhan2025} and \cite{FetcuOniciuc2022} is that they relied on a tensorial formula for a symmetric $(1,1)$ tensor field $S$ satisfying formally the Codazzi equation, that is $\nabla S$ is a symmetric $(1, 2)$ tensor field. To overcome this difficulty, we use a tensorial formula, introduced in \cite{FetcuLoubeauOniciuc2021}, which holds for any symmetric $(1,1)$ tensor field.
	
	More precisely, we consider a family of tensor fields, depending on a parameter $k$, whose vanishing divergence characterizes the biconservativity of the hypersurface. For a particular value of the parameter, the tensorial formula from \cite{FetcuLoubeauOniciuc2021} leads to a formula which is equivalent to the one obtained in \cite{AndronicKayhan2025}, and this formula provides a rigidity result only when $m \leq 10$. For the other values of the parameter, we obtain formulas which also hold when $m > 10$, provided that $|A|^2$ is bounded from above by a multiple of $f^2$. Then, we choose the parameter such that this upper bound is as large as possible.
	
	Our main result is the following. 
	\begin{customthm} {\ref{th:mainTheorem}}
		Let $\varphi : M^m \to N^{m+1} (c)$ be a compact non-minimal biconservative hypersurface with $m > 10$. If $\Riem^M \geq 0$ and the shape operator $A$ satisfies
		\begin{equation*}
			|A|^2 \leq \frac {m^2 (2m + 25)} {2 (m - 1)} f^2, \quad \text{on } M,
		\end{equation*}
		then $\nabla A = 0$. Moreover,
		\begin{enumerate} [label = \alph*)]
			\item if $c = -1$, then $\varphi (M)$ is a hypersphere $\mathbb S^m (r)$ of radius $r > 0$,
			
			\item if $c = 0$, then $\varphi (M)$ is a hypersphere $\mathbb S^m (r)$,
			
			\item if $c = 1$, then $\varphi (M)$ is either a small hypersphere $\mathbb S^m (r)$, $r \in (0, 1)$, or the extrinsic product $\mathbb S^{m_1} (r_1) \times \mathbb S^{m_2} (r_2)$, where 
			\begin{equation*}
				r_1^2 \in \left ( 0, \frac 1 {1 + x_2} \right ] \cup \left [ \frac 1 {1 + x_1}, 1 \right ),
			\end{equation*}
			$x_1$ and $x_2$, $x_1 < x_2$, are the solutions of the following second degree polynomial equation
			\begin{equation*}
				m_1 \left ( \frac {2m + 25} {2 (m - 1)} m_1 - 1 \right ) x^2 - \frac {2m + 25} {m - 1} m_1 m_2 x + m_2 \left ( \frac {2m + 25} {2 (m - 1)} m_2 - 1 \right ) = 0.
			\end{equation*}
		\end{enumerate}
	\end{customthm}
	If, in addition to the hypotheses of Theorem \ref{th:mainTheorem}, we assume $f \neq 0$ at any point of $M$, then the upper bound of $|A|^2$ can be replaced by the larger one
	\begin{equation*}
		\frac {m^2 (2m^2 + 5m - 34)} {2 (m - 1) (m - 10)} f^2
	\end{equation*}
	(see Theorem \ref{th:mainTheoremFneq0}). We also obtain local versions of these results, in which the compactness of $M$ is replaced by the condition that $|A|^2 + m^2 f^2 / 2$ is constant (see Theorems \ref{th:mainTheoremWithoutCompactness} and \ref{th:mainTheoremFneq0WithoutCompactness}).
	
	It is known that at any point where $\grad f \neq 0$, a biconservative hypersurface satisfies the J.H. Chen inequality (see \cite{JHChen1993})
	\begin{equation*}
		|A|^2 \geq \frac {m^2 (m + 8)} {4 (m - 1)} f^2.
	\end{equation*}
	We note that these upper bounds are compatible with the non-constancy of $f$, since they are larger than the lower bound of $|A|^2$ in J.H. Chen inequality. Therefore, the above pinching conditions do not imply, from the beginning, that $f$ is constant.
	
	The paper is organized as follows. In Section \ref{sec:preliminaries}, we recall some known results on biconservative hypersurfaces and the tensorial formula from \cite{FetcuLoubeauOniciuc2021}, and we derive a global expression of $\theta$ which appears in this formula. In Section \ref{sec:rigidity}, we introduce the tensors $T_k$, we use them in the tensorial formula and we obtain the main formula \eqref{eq:finalFormulaGenreal}. Then, we prove our rigidity results.
	
	\textbf{Conventions.}
	In this paper, we assume that all manifolds are connected and oriented. In general, the metric on an arbitrary manifold is denoted by $\langle \cdot, \cdot \rangle$ or not explicitly mentioned. The Levi-Civita connection of the Riemannian manifold $M$ is denoted by $\nabla$.
	
	The rough Laplacian acting on the set of all sections in the pull-back bundle $\varphi^{-1} (TN)$ is given by 
	$$
	\Delta^\varphi = - \trace \left ( \nabla ^\varphi \nabla ^\varphi - \nabla ^\varphi _\nabla \right ),
	$$
	where $\nabla ^\varphi$ is the connection in the pull-back bundle, and the curvature tensor field is 
	$$
	R (X, Y) Z = \nabla_X \nabla_Y Z - \nabla_Y \nabla_X Z - \nabla _{[X, Y]} Z.
	$$
	The mean curvature function of a hypersurface $\varphi : M^m \to N^{m+1}$ is denoted by $f = (\trace A) / m$, where $A = A_\eta$ is the shape operator of $M$ and $\eta$ is a unit section in the normal bundle.
	
	\section{Preliminaries} \label{sec:preliminaries}
	
	We recall that for a hypersurface $M^m$ in a space form $N^{m+1} (c)$ the Gauss equation is
	\begin{equation}\label{eq:GaussEquation}
		R(X, Y) Z = c \left ( \left \langle Y, Z \right \rangle X - \left \langle X, Z \right \rangle Y \right ) + \left \langle A(Y), Z \right \rangle A(X) - \left \langle A(X), Z \right \rangle A(Y)
	\end{equation}
	for any $X, Y, Z \in C(TM)$ and the Codazzi equation is
	\begin{equation} \label{eq:CodazziEquation}
		\left ( \nabla A\right ) (X, Y) = \left ( \nabla A \right ) (Y, X)
	\end{equation}
	for any $X, Y \in C(TM)$.
	
	Now we recall some well-known properties of the shape operator of hypersurfaces in space forms.
	\begin{lemma}\label{th:PropertiesNablaA}
		Let $\varphi : M^m \to  N^{m+1}(c)$ be a hypersurface in a space form. Then,
		\begin{enumerate}[label = \alph*)]
			\item $\left \langle (\nabla A) (\cdot, \cdot), (\cdot) \right \rangle$ is totally symmetric,
			\item $\left ( \Div A \right )^\sharp = \trace (\nabla A) (\cdot, \cdot) = m \grad f$,
			\item $\left ( \Div A^2 \right )^\sharp = m A(\grad f) + \grad |A|^2 / 2$. \label{item:divA2}
		\end{enumerate}
	\end{lemma}
	
	Consider a hypersurface $\varphi : M^m \to N^{m+1}(c)$ in a space form. The principal curvatures of $M$, denoted by $\{ \lambda_i \}_{i \in \overline{1, m}}$, are the eigenvalue functions of the shape operator. These functions are continuous on $M$ in general, but they need not be smooth. Let $M_A$ be the subset of points of $M$ around which the number of distinct principal curvatures is locally constant. It is a well-known fact that $M_A$ is an open and dense subset of $M$. On any connected component of $M_A$, which is itself open in $M$, the principal curvatures are smooth, and one can find a local orthonormal frame field $\{E_i\}_{i \in \overline{1, m}}$ satisfying $A(E_i) = \lambda_i E_i$ for all $i \in \overline{1, m}$. 
	
	Let $\varphi : M^m \to \left ( N^n, h \right )$ be a map, where $M$ is compact and the metric $h$ is fixed. On the set $\mathscr G$ of all metrics on $M$ we consider the functional
	\begin{equation*}
		\mathscr E_2 : \mathscr G \to \mathbb R, \quad \mathscr E_2 (g) = E_2 (\varphi).
	\end{equation*}
	In \cite{LoubeauMontaldoOniciuc2008} it was shown that its critical points are characterized by the vanishing of the stress-bienergy tensor field
	\begin{align*}
		S_2 (X, Y) =& \frac 1 2 |\tau (\varphi)|^2 \langle X, Y \rangle + \langle d \varphi, \nabla \tau (\varphi) \rangle \langle X, Y \rangle \\
		& - \langle d \varphi (X), \nabla _Y \tau (\varphi) \rangle - \langle d \varphi (Y), \nabla _X \tau (\varphi) \rangle,
	\end{align*}
	which was first introduced in \cite{Jiang1987} and satisfies
	\begin{equation*}
		\Div S_2 = - \langle \tau_2 (\varphi), d \varphi \rangle.
	\end{equation*}
	If $\varphi$ is an immersion, then $(\Div S_2)^\sharp = - \tau_2 (\varphi) ^\top$, where $\tau _2 (\varphi) ^\top$ denotes the tangent part of the bitension field. The biconservative submanifolds were defined in \cite{CaddeoMontaldoOniciucPiu2014} as those which satisfy $\Div S_2 = 0$. In the case of hypersurfaces, the stress-bienergy tensor field takes the form
	\begin{equation*}
		S_2 = 2mf A - \frac {m^2} 2 f^2 \id
	\end{equation*}
	and the biconservative equation becomes 
	\begin{equation}\label{eq:BiconservativityCharacterization}
		A(\grad f) = - \frac m 2 f \grad f.
	\end{equation}
	In \cite{AndronicKayhan2025}, the following characterization of biconservative hypersurfaces was obtained.
	\begin{lemma} \label{th:characterizationBiconservativityFirstPaper}
		A hypersurface $\varphi : M^m \to N^{m+1} (c)$ is biconservative if and only if $\Div \left ( f^2 A \right ) = 0$.
	\end{lemma}
	Now, we recall the well-known rigidity result obtained by Nomizu and Smyth in \cite{NomizuSmyth1969} which classifies compact hypersurfaces in space forms.
	\begin{theorem} \label{th:NomizuSmythTheorem}
		Let $\varphi : M^m \to N^{m+1} (c)$ be a compact hypersurface in a space form. If $M$ has non-negative sectional curvature, i.e. $\Riem ^M \geq 0$, and constant mean curvature, then $\nabla A = 0$ and $\varphi (M)$ is one of the following hypersurfaces
		\begin{enumerate}[label = \alph*)]
			\item the Euclidean hypersphere $\mathbb S^m (r)$ of radius $r > 0$, if $c \in \{-1, 0\}$, i.e. $N$ is either the hyperbolic space $\mathbb H^{m+1}$ or the Euclidean space $\mathbb R^{m+1}$;
			\item either the hypersphere $\mathbb S^m (r)$, $r \in (0, 1]$, or the extrinsic product $\mathbb S^{m_1} (r_1) \times \mathbb S^{m_2} (r_2)$, where $r_1^2 + r_2^2 = 1$, $m_1 + m_2 = m$, if $c = 1$, i.e. $N$ is the unit Euclidean sphere $\mathbb S^{m+1}$.
		\end{enumerate}
		Moreover, in the case a), as well as in the case b) when $\varphi (M)$ is a hypersphere or when $m_1 \geq 2$ and $m_2 \geq 2$, the immersion $\varphi$ has to be an embedding.
	\end{theorem}
	\begin{remark} \label{rem:localNomizuSmyth}
		In \cite{NomizuSmyth1969}, it was also shown that, in Theorem \ref{th:NomizuSmythTheorem}, the hypothesis of compactness can be replaced by the condition that $|A|^2$ is constant, and the conclusion $\nabla A = 0$ still holds.
	\end{remark}
	Biconservative hypersurfaces form a natural class in the geometry of submanifolds, as they allow this result to be generalized. Indeed, as shown in \cite{FetcuOniciuc2022}, if the CMC hypothesis in Theorem \ref{th:NomizuSmythTheorem} is replaced by biconservativity when $m \leq 10$, then the following result holds.
	\begin{theorem} \label{th:FetcuOniciuc2022}
		Let $\varphi : M^m \to N^{m+1} (c)$ be a compact non-minimal biconservative hypersurface. If $\Riem^M \geq 0$ and $m \leq 10$, then $\nabla A = 0$ and $\varphi (M)$ is one of the following hypersurfaces
		\begin{enumerate}[label = \alph*)]
			\item the Euclidean hypersphere $\mathbb S^m (r)$ of radius $r > 0$, if $c \in \{-1, 0\}$, i.e. $N$ is either the hyperbolic space $\mathbb H^{m+1}$ or the Euclidean space $\mathbb R^{m+1}$;
			\item either the small hypersphere $\mathbb S^m (r)$, $r \in (0, 1)$, or the extrinsic product $\mathbb S^{m_1} (r_1) \times \mathbb S^{m_2} (r_2)$, where $r_1^2 + r_2^2 = 1$, $m_1 + m_2 = m$ and $r_1 \neq \sqrt {m_1 / m}$, if $c = 1$, i.e. $N$ is the unit Euclidean sphere $\mathbb S^{m+1}$.
		\end{enumerate}
	\end{theorem}
	
	Now, we recall some well-known properties of biconservative hypersurfaces in space forms (see, for example, \cite{BalmusMontaldoOniciuc2008}, \cite{BibiLoubeauOniciuc2021} and \cite{Fu2015Nachr}).
	
	\begin{lemma} \label{th:LemmaBiconservativeHypersurfacesConnectionForms}
		Let $\varphi : M^m \to N^{m+1} (c)$ be a biconservative hypersurface in a space form. Assume that $\grad f \neq 0$ at any point of $M$ and the number of distinct principal curvatures is constant. Then, $\lambda_1 = - m f / 2$ is a principal curvature of $A$, globally defined, with multiplicity $m_1 = 1$ and $E_1 = \grad f / |\grad f|$. Moreover, $E_a (f) = 0$, for any $a \in \overline {2, m}$ and $\nabla _{E_1} E_1 = 0$.
	\end{lemma}	
	
	In Proposition 4.1 of \cite{FetcuLoubeauOniciuc2021}, the following formula was obtained for a symmetric $(1, 1)$ tensor $S$ defined on a manifold $M$
	\begin{equation} \label{eq:mainTensorialFormulaInitial}
		\Div \theta = - |\nabla S|^2 + |\Div S|^2 + \frac 1 2 |W|^2 + \langle T, S \rangle,
	\end{equation}
	where $T (X) = - \trace (RS) (\cdot, X, \cdot)$, $W (X, Y) = (\nabla S) (X, Y) - (\nabla S) (Y, X)$, for any $X, Y \in C(TM)$ and $\theta$ is defined locally by $\theta^i = \left ( g^{ik} g^{j\ell} - g^{jk} g^{i \ell} \right ) \left ( \nabla _j S^a_\ell \right ) S_{ak}$.
	
	Let $\{ e_i \}_{i \in \overline {1, m}}$ be an orthonormal basis in $T_p M$, $p \in M$, such that $S (e_i) = \mu_i e_i$, for any $i \in \overline {1, m}$. In this orthonormal basis, the last term of \eqref{eq:mainTensorialFormulaInitial} becomes
	\begin{equation*}
		\langle T, S \rangle = - \frac 1 2 \sum _{i, j = 1} ^m (\mu_i - \mu_j)^2 R_{ijij},
	\end{equation*}
	where $R_{ijij} = \langle R(e_i, e_j) e_j, e_i \rangle$. Thus, we obtain
	\begin{equation} \label{eq:mainTensorialFormula}
		\Div \theta = - |\nabla S|^2 + |\Div S|^2 + \frac 1 2 |W|^2 - \frac 1 2 \sum _{i, j = 1} ^m (\mu_i - \mu_j)^2 R_{ijij}.
	\end{equation}
	We derive now a global expression of $\theta$ in terms of $S$, which will be useful later. Using the fact that $\nabla g = 0$, we obtain 
	\begin{align}
		\theta^i&=\bigl(g^{ik}g^{j\ell}-g^{jk}g^{i\ell}\bigr)(\nabla_j S^a_\ell)\,S_{ak} \label{eq:expressionThetaI}\\
				&=g^{j\ell}(\nabla_j S^a_\ell)\,S^i_a-g^{i\ell}(\nabla_j S^a_\ell)\,S^j_a \notag\\
				&=g^{j\ell}(\nabla_j S^a_\ell)\,S^i_a-(\nabla_j S^{ia})\,S^j_a \notag\\
				&=g^{j\ell}(\nabla_j S^a_\ell)\,S^i_a-g^{a\ell}(\nabla_j S^i_\ell)\,S^j_a \notag\\
				&=g^{j\ell}(\nabla_j S^a_\ell)\,S^i_a-S^{\ell j}(\nabla_j S^i_\ell). \notag
	\end{align}
	We note that 
	\begin{equation*}
		\left ( \Div S \right )^\sharp = \trace \nabla S = g^{ij} \left ( \nabla _i S^k_j \right ) \partial_k
	\end{equation*}
	and
	\begin{equation*}
		S \left ( \left ( \Div S \right )^\sharp \right ) = g^{ij} S^\ell_k \left ( \nabla _i S^k_j \right ) \partial_k.
	\end{equation*}
	We have
	\begin{align*}
		\left ( \Div S^2 \right )^\sharp =& \trace \nabla S^2 \\
				 						 =& g^{ij} \left ( \nabla _i \left ( S^a_j S^k_a \right ) \right ) \partial_k \\
				 						 =& g^{ij} \left ( S^k_a \left ( \nabla_i S^a_j \right ) + S^a_j \left ( \nabla _i S^k_a \right ) \right ) \partial_k \\
				 						 =& S \left ( \left ( \Div S \right ) ^\sharp \right ) + S^{ia} \left ( \nabla _i S^k_a \right ) \partial_k,
	\end{align*}
	which implies 
	\begin{equation*}
		S^{ia} \left ( \nabla _i S^k_a \right ) \partial_k = \left ( \Div S^2 \right )^\sharp - S \left ( \left ( \Div S \right ) ^\sharp \right ).
	\end{equation*}
	Therefore, from \eqref{eq:expressionThetaI}, we obtain
	\begin{equation} \label{eq:formulaTheta}
		\theta = 2 S \left ( \left ( \Div S \right ) ^\sharp \right ) - \left ( \Div S^2 \right )^\sharp.
	\end{equation}
	
	\section{Rigidity of biconservative hypersurfaces} \label{sec:rigidity}
	
	In this section we use \eqref{eq:mainTensorialFormula} to derive new rigidity results for compact biconservative hypersurfaces in space forms. To make full use of \eqref{eq:mainTensorialFormula}, we introduce a new family of tensors whose vanishing divergence characterizes the biconservativity of a hypersurface.
	
	\begin{proposition} \label{th:newTensorsWhichCharacterizeBiconservativity}
		Let $k$ be a positive integer. A hypersurface $\varphi : M^m \to N^{m+1} (c)$ is biconservative if and only if the symmetric $(1, 1)$ tensor 
		\begin{equation*}
			T_k = f^k A + \frac {m (k - 2)} {2 (k + 1)} f^{k+1} \id, 
		\end{equation*}
		is divergence-free.
	\end{proposition}
	
	\begin{proof}
		First, we compute the divergence of $T_k$. For this we consider $\left \{ X_i \right \}_{i \in \overline {1, m}}$ a local orthonormal frame field. We have 
		\begin{align*}
			\left ( \Div T_k \right )^\sharp =& \left ( \Div \left ( f^k A + \frac {m (k - 2)} {2 (k + 1)} f^{k+1} \id \right ) \right )^\sharp\\
					 =& \sum _{i = 1} ^m \left ( \nabla \left ( f^k A + \frac {m (k - 2)} {2 (k + 1)} f^{k+1} \id \right ) \right ) (X_i, X_i) \\
					 =& \sum _{i = 1} ^m \left ( \nabla _{X_i} \left ( f^k A X_i \right ) - f^k A \left ( \nabla _{X_i} X_i \right ) + \frac {m (k - 2)} {2 (k + 1)} \left ( \nabla _{X_i} \left ( f^{k+1} X_i \right ) - f^{k + 1} \nabla _{X_i} X_i \right ) \right ) \\
					 =& \sum _{i = 1} ^m \left ( X_i \left ( f^k \right ) AX_i + f^k \nabla _{X_i} (AX_i) - f^k A \left ( \nabla _{X_i} X_i \right ) + \frac {m (k - 2)} {2 (k + 1)} X_i \left ( f^{k+1} \right ) X_i \right ) \\
					 =& A \left ( \grad f^k \right ) + f^k \trace (\nabla A) + \frac {m (k - 2)} {2 (k + 1)} \grad f^{k+1} \\
					 =& k f^{k - 1} A (\grad f) + m f^k \grad f + \frac {m (k - 2)} 2 f^k \grad f \\
					 =& k f^{k - 1} A (\grad f) + \frac {m k} 2 f^k \grad f.
		\end{align*}
		Thus, 
		\begin{equation*}
			\left ( \Div T_k \right ) ^\sharp = k f^{k - 1} \left ( A (\grad f) + \frac m 2 f \grad f \right ).
		\end{equation*}
		It is known that 
		\begin{equation*}
			\left ( \Div S_2 \right )^\sharp = 2m \left ( A (\grad f) + \frac m 2 f \grad f \right ).
		\end{equation*}
		Obviously, if $\varphi$ is biconservative, then $\Div T_k = 0$.
		
		Suppose that $\Div T_k = 0$. Then, at any point,
		\begin{equation*}
			f = 0 \quad \text{or} \quad A (\grad f) + \frac m 2 f \grad f = 0.
		\end{equation*}
		If $A (\grad f) + m f / 2 \grad f = 0$, then $\Div S_2 = 0$, i.e. $\varphi$ is biconservative.
		
		Assume by way of contradiction that $A (\grad f) + m f / 2 \grad f \neq 0$ on an open subset $U$ of $M$. Then, $f = 0$ on $U$ and, consequently, $\grad f = 0$ on $U$. This implies that $A (\grad f) + m f / 2 \grad f = 0$, contradiction.
	\end{proof}
	
	\begin{remark}
		Note that 
		\begin{equation*}
			T_1 = \frac 1 {2m} \left ( 2mfA - \frac {m^2} 2 f^2 \id \right ) = \frac 1 {2m} S_2 \qquad \text{and} \qquad T_2 = f^2 A.
		\end{equation*}
		Therefore, this family of tensors is a good generalization of the stress-bienergy tensor and the tensor from Lemma \ref{th:characterizationBiconservativityFirstPaper}.
	\end{remark}
	
	\begin{remark}
		Since the mean curvature function $f$ can vanish at some points, we cannot consider $k = 0$ in Proposition \ref{th:newTensorsWhichCharacterizeBiconservativity}. Instead, we can define 
		$$
		T_0 = A - mf \id,
		$$
		which is divergence-free independent of the biconservativity condition. Then, we can consider the tensor family $\{ T_k \}_{k \geq 1} \cup \{ T_0 \}$ and study what happens in the case $k = 0$.
	\end{remark}
	
	We now give some general properties of the tensors $T_k$, which hold for any hypersurface.
	\begin{lemma} \label{th:nablaTk}
		Let $\varphi : M^m \to N^{m+1} (c)$ be a hypersurface in a space form and $k$ a positive integer. Then, for any $X, Y \in C(TM)$, we have
		\begin{enumerate}[label = \alph*)]
			\item \label{item:Tka} the covariant derivative of $T_k$ is given by
			\begin{equation*}
				(\nabla T_k) (X, Y) = f^k (\nabla A) (X, Y) + X \left ( f^k \right ) A (Y) + \frac {m (k - 2)} {2 (k + 1)} X \left ( f^{k+1} \right ) Y,
			\end{equation*}
			\item \label{item:Tkb} the tensor field $\nabla T_k$ satisfies
			\begin{align*}
				(\nabla T_k) (X, Y) =& (\nabla T_k) (Y, X) + A \left ( X \left ( f^k \right ) Y - Y \left ( f^k \right ) X \right ) \\
									 & + \frac {m (k - 2)} {2 (k + 1)} \left ( X \left ( f^{k+1} \right ) Y - Y \left ( f^{k+1} \right ) X \right ).
			\end{align*}
		\end{enumerate}
	\end{lemma}
	
	\begin{proof}
		For any $X, Y \in C(TM)$, we have
		\begin{align*}
			\left ( \nabla T_k \right ) (X, Y) =& \nabla _X (T_k (Y)) - T_k \left ( \nabla _X Y \right ) \\
											   =& \nabla _X \left ( f^k A (Y) + \frac {m (k - 2)} {2 (k + 1)} f^{k+1} Y \right ) - f^k A \left ( \nabla _X Y \right ) - \frac {m (k - 2)} {2 (k + 1)} f^{k+1} \nabla _X Y \\
											   =& X \left ( f^k \right ) A (Y) + f^k \nabla _X (A(Y)) + \frac {m (k - 2)} {2 (k + 1)} X \left ( f^{k+1} \right ) Y + \frac {m (k - 2)} {2 (k + 1)} f^{k+1} \nabla _X Y \\
											    & - f^k A \left ( \nabla _X Y \right ) - \frac {m (k - 2)} {2 (k + 1)} f^{k+1} \nabla _X Y.
		\end{align*}
		Since $\nabla_X (A(Y)) - A \left ( \nabla _X Y \right ) = (\nabla A) (X, Y)$, we obtain \ref{item:Tka}. Relation \ref{item:Tkb} follows from \ref{item:Tka} and the Codazzi equation \eqref{eq:CodazziEquation}.
	\end{proof}
	
	In the following, we consider $\varphi : M^m \to N^{m+1} (c)$ a biconservative hypersurface in a space form and use the tensors $T_k$ in formula \eqref{eq:mainTensorialFormula}. We denote by $W_k$ and $\theta_k$ the tensor field $W$ and the vector field $\theta$ from \eqref{eq:mainTensorialFormulaInitial} corresponding to $S = T_k$.
	
	First, we want to obtain the relation between the eigenvalue functions $\mu_i$ of $T_k$ and the principal curvatures $\lambda_i$. For this, we consider an orthonormal basis $\{e_i\}_{i \in \overline {1, m}}$ in $T_p M$, $p \in M$, such that $A (e_i) = \lambda_i e_i$. It is easy to see that $\{ e_i \}_{i \in \overline {1, m}}$ also diagonalizes $T_k$. Thus, 
	\begin{equation*}
		\mu_i = \langle T_k (e_i), e_i \rangle = \left \langle f^k A (e_i) + \frac {m (k - 2)} {2 (k + 1)} f^{k+1} e_i, e_i \right \rangle, 
	\end{equation*}
	which implies that 
	\begin{equation*}
		\mu_i = f^k \left ( \lambda_i + \frac {m (k - 2)} {2 (k + 1)} f \right ).
	\end{equation*}
	Therefore, 
	\begin{equation} \label{eq:lastTermOfMainFormulaTk}
		- \frac 1 2 \sum _{i, j = 1} ^m (\mu_i - \mu_j)^2 R_{ijij} = - \frac 1 2 f^{2k} \sum _{i, j = 1} ^m (\lambda_i - \lambda_j)^2 R_{ijij}.
	\end{equation}
	Using item \ref{item:Tkb} of Lemma \ref{th:nablaTk}, we obtain 
	\begin{align*}
		W_k (X, Y) =& (\nabla T_k) (X, Y) - (\nabla T_k) (Y, X) \\
				   =& X \left ( f^k \right ) A(Y) - Y \left ( f^k \right ) A (X) + \frac {m (k - 2)} {2 (k + 1)} \left ( X \left ( f^{k+1} \right ) Y - Y \left ( f^{k+1} \right ) X \right ).
	\end{align*}
	Now, we consider a local orthonormal frame field $\left \{ X_i \right \} _{i \in \overline {1, m}}$ and compute
	\begin{align*}
		\frac 1 2 |W_k|^2 =& \frac 1 2 \sum _{i, j = 1} ^m |W_k (X_i, X_j)|^2 \\
			  			=& \frac 1 2 \sum _{i, j = 1} ^m \left | X_i \left ( f^k \right ) A(X_j) - X_j \left ( f^k \right ) A(X_i) + \frac {m (k - 2)} {2 (k + 1)} \left ( X_i \left ( f^{k + 1} \right ) X_j - X_j \left ( f^{k + 1} \right ) X_i \right ) \right |^2 \\
			  			=& \frac 1 2 \sum _{i, j = 1} ^m \left [ \left ( X_i \left ( f^k \right ) \right )^2 |A (X_j)|^2 - 2 X_i \left ( f^k \right ) X_j \left ( f^k \right ) \langle A (X_j), A(X_i) \rangle \vphantom{\frac {2 m (k - 2)} {2 (k + 1)}} \right . \\
					     & \qquad \quad + \frac {2 m (k - 2)} {2 (k + 1)} X_i \left ( f^k \right ) X_i \left ( f^{k + 1} \right ) \langle A (X_j), X_j \rangle \\
					     & \qquad \quad - \frac {2 m (k - 2)} {2 (k + 1)} X_i \left ( f^k \right ) X_j \left ( f^{k+1} \right ) \langle A (X_j), X_i \rangle \\
					     & \qquad \quad + \left ( X_j \left ( f^k \right ) \right )^2 \left | A (X_i) \right |^2 - \frac {2 m (k - 2)} {2 (k + 1)} X_j \left ( f^k \right ) X_i \left ( f^{k + 1} \right ) \langle A(X_i), X_j \rangle \\
					     & \qquad \quad + \frac {2 m (k - 2)} {2 (k + 1)} X_j \left ( f^k \right ) X_j \left ( f^{k + 1} \right ) \langle A (X_i), X_i \rangle \\
					     & \qquad \quad + \frac {m^2 (k - 2)^2} {4 (k + 1)^2} \left ( X_i \left ( f^{k+1} \right ) \right )^2 |X_j|^2 - \frac {2 m^2 (k - 2)^2} {4 (k + 1)^2} X_i \left ( f^{k+1} \right ) X_j \left ( f^{k + 1} \right ) \langle X_i, X_j \rangle \\
					     & \left . \qquad \quad + \frac {m^2 (k - 2)^2} {4 (k + 1)^2} \left ( X_j \left ( f ^{k + 1} \right ) \right )^2 |X_i|^2 \right ] 
	\end{align*}
	Thus, we obtain
	\begin{align}
		\frac 1 2 |W_k|^2 =& |A|^2 \left | \grad f^k \right |^2 - \left | A \left ( \grad f^k \right ) \right |^2 + \frac {m^2 (k - 2)} {k + 1} f \left \langle \grad f^k, \grad f^{k+1} \right \rangle \label{eq:NormWBeforeDerivation} \\
		& - \frac {m (k - 2)} {k + 1} \left \langle A \left ( \grad f^{k+1} \right ), \grad f^k \right \rangle + \frac {m^2 (m - 1) (k - 2)^2} {4 (k + 1)^2} \left | \grad f^{k+1} \right |^2. \notag
	\end{align}
	When $k \geq 1$, we use the biconservativity of $M$ in \eqref{eq:NormWBeforeDerivation} and obtain
	\begin{align*}
		\frac 1 2 |W_k|^2 =& k^2 f^{2 (k - 1)} |A|^2 |\grad f|^2 - \frac {m^2 k^2} 4 f^{2k} |\grad f|^2 \\
						 & + m^2 k (k - 2) f^{2k} |\grad f|^2 + \frac {m^2 k (k - 2)} 2 f^{2k} |\grad f|^2 \\
						 & + \frac {m^2 (m - 1) (k - 2)^2} 4 f^{2k} |\grad f|^2 \\
						=& k^2 f^{2 (k - 1)} |A|^2 |\grad f|^2 \\
						 & + \left ( - \frac {m^2 k^2} 4 + m^2 k (k - 2) + \frac {m^2 k (k - 2)} 2 + \frac {m^2 (m - 1) (k - 2)^2} 4 \right ) f^{2k} |\grad f|^2,
	\end{align*}
	which is equivalent to
	\begin{equation} \label{eq:NormaWGenral}
		\frac 1 2 |W_k|^2 = k^2 f^{2 (k - 1)} |A|^2 |\grad f|^2 + \frac {m^2 \left ( (m + 4) k^2 - 4(m + 2) k + 4 (m - 1) \right )} 4 f^{2k} |\grad f|^2.
	\end{equation}
	
	Now, we want to find a lower bound for $|\nabla T_k|^2$ in terms of $|\grad f|^2$. For this, we consider the set $M_A$ of all points of $M$ at which the number of distinct principal curvatures is locally constant. Let $U$ be a connected component of $M_A$ and $p \in U$.
	
	Assume that $(\grad f) (p) \neq 0$. Since $U$ is open, we deduce that there is an open subset $V \subset U$, $p \in V$, such that 
	\begin{equation*}
		\grad f \neq 0, \quad \text{at any point of } V.
	\end{equation*}
	Since $p \in U$, eventually restricting $V$, there exists a local orthonormal frame field $\{ E_i \}_{i \in \overline {1, m}}$ defined on $V$ which diagonalizes the shape operator $A$. Since $\grad f \neq 0$ at any point of $V$ and the biconservative condition \eqref{eq:BiconservativityCharacterization}, we can assume that 
	\begin{equation*}
		E_1 = \frac {\grad f} {|\grad f|} \quad \text{and} \quad \lambda_1 = - \frac {m f} 2.
	\end{equation*}
	From Lemma \ref{th:LemmaBiconservativeHypersurfacesConnectionForms}, the multiplicity of $\lambda_1$ is $m_1 = 1$.
	
	We have 
	\begin{align*}
		|\nabla T_k|^2 =& \sum _{i, j = 1} ^m |(\nabla T_k) (E_i, E_j)|^2 = \sum _{i, j, \ell = 1} ^m \left | \langle (\nabla T_k) (E_i, E_j), E_\ell \rangle E_\ell \right |^2 \\
					   =& \sum _{i, j, \ell = 1} ^m \langle (\nabla T_k) (E_i, E_j), E_\ell \rangle ^2 \\
					\geq& \langle (\nabla T_k) (E_1, E_1), E_1 \rangle ^2 + \sum _{i = 2} ^m \langle (\nabla T_k) (E_i, E_i), E_1 \rangle ^2 \\
						& + \sum _{i = 2} ^m \langle (\nabla T_k) (E_i, E_1), E_i \rangle ^2 + \sum_{i = 2} ^m \langle (\nabla T_k) (E_1, E_i), E_i \rangle ^2 \\
					   =& \langle (\nabla T_k) (E_1, E_1), E_1 \rangle ^2 + 2 \sum _{i = 2} ^m \langle (\nabla T_k) (E_i, E_i), E_1 \rangle ^2 \\
					    &  + \sum_{i = 2} ^m \langle (\nabla T_k) (E_1, E_i), E_i \rangle ^2 \\
					\geq& \langle (\nabla T_k) (E_1, E_1), E_1 \rangle ^2 + \frac 2 {m - 1} \left ( \sum _{i = 2} ^m \langle (\nabla T_k) (E_i, E_i), E_1 \rangle \right )^2 \\
					    & + \frac 1 {m - 1} \left ( \sum _{i = 2} ^m \langle (\nabla T_k) (E_1, E_i), E_i \rangle \right )^2.
	\end{align*}
	In the above computation, we used the fact that, since $T_k$ is symmetric, $\langle (\nabla T_k) (X, Y), Z \rangle = \langle (\nabla T_k) (X, Z), Y \rangle$, for any $X, Y, Z \in C(TM)$, and thus $\langle (\nabla T_k) (E_i, E_1), E_i \rangle = \langle (\nabla T_k) (E_i, E_i), E_1 \rangle$.

	In the following, we compute each term of the right-hand side separately. First, using Lemma \ref{th:LemmaBiconservativeHypersurfacesConnectionForms}, we obtain
	\begin{align*}
		\left \langle (\nabla T_k) (E_1, E_1), E_1 \right \rangle =& \left \langle \nabla _{E_1} (T_k (E_1)) - T_k \left ( \nabla _{E_1} E_1 \right ), E_1 \right \rangle \\
								=& \left \langle \nabla _{E_1} \left ( f^k A (E_1) + \frac {m (k - 2)} {2 (k + 1)} f^{k+1} E_1 \right ), E_1 \right \rangle \\
								=& \left \langle \nabla _{E_1} \left ( - \frac m 2 f^{k+1} E_1 + \frac {m (k - 2)} {2 (k + 1)} f^{k+1} E_1 \right ), E_1 \right \rangle \\
								=& \left \langle \nabla _{E_1} \left ( - \frac {3m} {2 (k + 1)} f^{k+1} E_1 \right ), E_1 \right \rangle \\
								=& - \frac {3m} {2 (k + 1)} \left \langle E_1 \left ( f^{k+1} \right ) E_1 + f^{k+1} \nabla_{E_1} E_1, E_1 \right \rangle .
	\end{align*}
	Using Lemma \ref{th:LemmaBiconservativeHypersurfacesConnectionForms}, this relation is equivalent to
	\begin{equation} \label{eq:NablaTk111BeforeDerivation}
		\left \langle (\nabla T_k) (E_1, E_1), E_1 \right \rangle = - \frac {3m} {2 (k + 1)} E_1 \left (f^{k+1} \right ).
	\end{equation}
	Using Proposition \ref{th:newTensorsWhichCharacterizeBiconservativity} and \eqref{eq:NablaTk111BeforeDerivation}, we have 
	\begin{align*}
		\sum _{i = 2} ^m \langle (\nabla T_k) (E_i, E_i), E_1 \rangle =& \sum _{i = 1} ^m \langle (\nabla T_k) (E_i, E_i), E_1 \rangle - \langle (\nabla T_k) (E_1, E_1), E_1 \rangle \\
							=& \langle \left ( \Div T_k \right )^\sharp, E_1 \rangle - \langle (\nabla T_k) (E_1, E_1), E_1 \rangle \\
							=& - \langle (\nabla T_k) (E_1, E_1), E_1 \rangle,
	\end{align*}
	which is equivalent to
	\begin{equation} \label{eq:NablaTkii1BeforeDerivation}
		\sum _{i = 2} ^m \langle (\nabla T_k) (E_i, E_i), E_1 \rangle = \frac {3m} {2 (k + 1)} E_1 \left ( f^{k+1} \right ).
	\end{equation}
	In order to compute the last term, we use item \ref{item:Tkb} of Lemma \ref{th:nablaTk} and obtain 
	\begin{align*}
		\sum _{i = 2} ^m \left \langle (\nabla T_k) (E_1, E_i), E_i \right \rangle =& \sum _{i = 2} ^m \left \langle (\nabla T_k) (E_i, E_1) + A \left ( E_1 \left ( f^k \right ) E_i - E_i \left ( f^k \right ) E_1 \right ) \vphantom{\frac {m (k - 2)} {2(k + 1)}} \right . \\
							 & \left . \qquad + \frac {m (k - 2)} {2(k + 1)} \left ( E_1 \left ( f^{k+1} \right ) E_i - E_i \left ( f^{k+1} \right ) E_1 \right ), E_i \right \rangle \\
							=& \sum _{i = 2} ^m \left ( \langle (\nabla T_k) (E_i, E_1), E_i \rangle + E_1 \left ( f^k \right ) \langle A(E_i), E_i \rangle + \frac {m (k - 2)} {2 (k + 1)} E_1 \left ( f^{k+1} \right ) \right ) \\
							=& \sum_{i = 1} ^m \left ( \langle (\nabla T_k) (E_i, E_i), E_1 \rangle + E_1 \left ( f^k \right ) \langle A(E_i), E_i \rangle \right )\\
							 & - \langle (\nabla T_k) (E_1, E_1), E_1 \rangle - E_1 \left ( f^k \right ) \langle A (E_1), E_1 \rangle \\
							 & + \frac {m (m - 1) (k - 2)} {2 (k + 1)} E_1 \left ( f^{k + 1} \right ) \\
							=& \langle \left ( \Div T_k \right )^\sharp, E_1 \rangle + E_1 \left ( f^k \right ) \trace A - \langle (\nabla T_k) (E_1, E_1), E_1 \rangle + \frac {mf} 2 E_1 \left ( f^k \right ) \\
							 & + \frac {m (m - 1) (k - 2)} {2 (k + 1)} E_1 \left ( f^{k+1} \right ),
	\end{align*}
	which, since $\Div T_k = 0$, is equivalent to
	\begin{equation} \label{eq:NablaTk1iiBeforeDerivation}
		\sum _{i = 2} ^m \left \langle (\nabla T_k) (E_1, E_i), E_i \right \rangle = \frac {3m} 2 f E_1 \left ( f^k \right ) + \frac {m \left ( (m - 1) k - 2 m + 5 \right )} {2 (k + 1)} E_1 \left ( f^{k+1} \right ).
	\end{equation}
	Using \eqref{eq:NablaTk111BeforeDerivation}, \eqref{eq:NablaTkii1BeforeDerivation} and \eqref{eq:NablaTk1iiBeforeDerivation}, we obtain
	\begin{align}
		|\nabla T_k|^2 \geq & \frac {m^2 \left ( 9m + 9 + \left ( (m - 1) k - 2m + 5 \right )^2 \right )} {4 (m - 1) (k + 1)^2} \left ( E_1 \left ( f^{k+1} \right ) \right )^2 + \frac {9m^2} {4 (m - 1)} f^2 \left ( E_1 \left ( f^k \right ) \right )^2 \label{eq:lowerBoundNormNablaTkBeforeDerivation} \\
						   & + \frac {3m^2 \left ( (m - 1)k - 2m + 5 \right )} {2 (m - 1) (k + 1)} f E_1 \left ( f^k \right ) E_1 \left ( f^{k+1} \right ). \notag
	\end{align}
	When $k \geq 1$, taking into account the definition of $E_1$, from \eqref{eq:lowerBoundNormNablaTkBeforeDerivation}, we obtain that at the points where $\grad f \neq 0$, we have
	\begin{equation} \label{eq:lowerBoundNablaTkGeneral}
		|\nabla T_k|^2 \geq \frac {m^2 \left ( 9m + 9 + \left ( (m + 2) k - 2m + 5 \right )^2 \right )} {4 (m - 1)} f^{2k} |\grad f|^2.
	\end{equation}
	Obviously, relation \eqref{eq:lowerBoundNablaTkGeneral} also holds at the points where $\grad f = 0$. Therefore, relation \eqref{eq:lowerBoundNablaTkGeneral} holds on $M_A$, and consequently on the whole $M$.
	
	Using \eqref{eq:lastTermOfMainFormulaTk}, \eqref{eq:NormaWGenral} and \eqref{eq:lowerBoundNablaTkGeneral} in \eqref{eq:mainTensorialFormula} and taking into account that $\Div T_k = 0$, we obtain the main formula of this paper
	\begin{align}
		\Div \theta_k \leq& k^2 f^{2 (k - 1)} |\grad f|^2 \left ( |A|^2 - \left ( - \frac {3m^2 (m - 10)} {4 (m - 1)} \frac 1 {k^2} + \frac {3m^2 (m + 2)} {2 (m - 1)} \frac 1 k + \frac {m^2 (m + 8)} {4 (m - 1)} \right ) f^2 \right ) \label{eq:finalFormulaGenreal} \\
						  & - \frac 1 2 f^{2k} \sum _{i, j = 1} ^m (\lambda_i - \lambda_j)^2 R_{ijij}, \notag 
	\end{align}
	for any $k \geq 1$.
	\begin{remark}
		In the definition of the tensors $T_k$, the coefficient of $f^{k+1} \id$ was chosen such that $T_k$ is divergence-free. One can also consider, for any real number $\alpha_k$, the symmetric $(1, 1)$ tensor field
		\begin{equation*}
			T_{k, \alpha_k} = f^k A + \alpha_k f^{k+1} \id,
		\end{equation*}
		which, in general, is not divergence-free. Following the same steps as above and taking into account that, in this case, the term $|\Div T_{k, \alpha_k}|^2$ in \eqref{eq:mainTensorialFormula} does not vanish, by straightforward computations, one obtains for $\Div \theta_{k, \alpha_k}$ the same upper bound as in \eqref{eq:finalFormulaGenreal}, for any $\alpha_k \in \mathbb R$. Therefore, the choice $\alpha_k = m (k - 2) / (2 (k + 1))$ is not a restriction, and the upper bounds of $|A|^2$ obtained below cannot be improved by using the tensors $T_{k, \alpha_k}$.
	\end{remark}
	\begin{remark}
		To study formula \eqref{eq:mainTensorialFormula} for the tensor $T_0$, we do not need to repeat the computations from scratch. The expression for $|W_0|^2 / 2$ follows from \eqref{eq:NormWBeforeDerivation} by formally replacing $k = 0$, with the convention that $f^0 = 1$, including at points where $f = 0$. Therefore,
		\begin{equation} \label{eq:NormWK=0}
			\frac 1 2 |W_0|^2 = m^2 (m - 1) |\grad f|^2.	
		\end{equation}
		Again, we can obtain the lower bound of $|\nabla T_0|^2$ by	formally replacing $k = 0$ in \eqref{eq:lowerBoundNormNablaTkBeforeDerivation}. Thus, taking into account the definition of $E_1$, we obtain that, at the points where $\grad f \neq 0$, we have
		\begin{equation} \label{eq:lowerBoundNormNablaT0Final}
			\left | \nabla T_0 \right |^2 \geq \frac {m^2 \left ( 4m^2 - 11m + 34 \right )} {4 (m - 1)} |\grad f|^2.
		\end{equation}
		Obviously, relation \eqref{eq:lowerBoundNormNablaT0Final} also holds at the points where $\grad f = 0$. Therefore, relation \eqref{eq:lowerBoundNormNablaT0Final} holds on $M_A$ and consequently on the whole $M$. Now, replacing \eqref{eq:NormWK=0} and \eqref{eq:lowerBoundNormNablaT0Final} in \eqref{eq:mainTensorialFormula} and taking into account that $\Div T_0 = 0$, we obtain
		\begin{equation} \label{eq:finalFormulaT0}
			\Div \theta_0 \leq \frac {3m^2 (m - 10)} {4 (m - 1)} |\grad f|^2 - \frac 1 2 \sum _{i, j = 1} ^m (\lambda_i - \lambda_j)^2 R_{ijij}.
		\end{equation}
		Using \eqref{eq:formulaTheta}, it can be shown that \eqref{eq:finalFormulaT0} is equivalent to (2.8) from \cite{AndronicKayhan2025}, and therefore it does not provide additional information in our study.
	\end{remark}
	If $M$ is compact, then, after integrating, the left-hand side of \eqref{eq:finalFormulaGenreal} vanishes. In the following, we derive some conditions which force the right-hand side of \eqref{eq:finalFormulaGenreal} to be non-positive. Since the case when $m \leq 10$ was already studied in \cite{FetcuOniciuc2022}, we suppose that $m > 10$. If we assume that $\Riem^M \geq 0$, then the right-hand side of \eqref{eq:finalFormulaGenreal} is non-positive when, on $M$, we have
	\begin{equation} \label{eq:upperBoundNormAInitial}
		|A|^2 \leq \left ( - \frac {3m^2 (m - 10)} {4 (m - 1)} \frac 1 {k^2} + \frac {3m^2 (m + 2)} {2 (m - 1)} \frac 1 k + \frac {m^2 (m + 8)} {4 (m - 1)} \right ) f^2.
	\end{equation}
	Now, we find the maximum upper bound for $|A|^2$. To do this, we denote 
	\begin{equation*}
		b_k = -3 (m - 10) \frac 1 {k^2} + 6 (m + 2) \frac 1 k + m + 8,
	\end{equation*}
	and rewrite \eqref{eq:upperBoundNormAInitial} as 
	\begin{equation*}
		|A|^2 \leq \frac {m^2} {4 (m - 1)} b_k f^2.
	\end{equation*}
	Since the leading coefficient of $b_k$ is negative, the maximum value of $b_k$ is 
	\begin{equation*}
		b_{k_0} = \frac {2\left ( 2m^2 + 5m - 34 \right )} {m - 10},
	\end{equation*}
	where $k_0 = (m - 10) / (m + 2) < 1$. Since we study the case when $k$ is integer, we deduce that the maximum upper bound of $|A|^2$ is attained when $k = 1$ and, in this case, we have
	\begin{equation} \label{eq:upperBoundNormAForK=1}
		|A|^2 \leq \frac {m^2 (2m + 25)} {2 (m - 1)} f^2.
	\end{equation}
	Based on this, we have the following result:
	\begin{theorem} \label{th:mainTheorem}
		Let $\varphi : M^m \to N^{m+1} (c)$ be a compact non-minimal biconservative hypersurface with $m > 10$. If $\Riem^M \geq 0$ and \eqref{eq:upperBoundNormAForK=1} holds on $M$, then $\nabla A = 0$. Moreover,
		\begin{enumerate} [label = \alph*)]
			\item if $c = -1$, then $\varphi (M)$ is a hypersphere $\mathbb S^m (r)$ of radius $r > 0$,
			
			\item if $c = 0$, then $\varphi (M)$ is a hypersphere $\mathbb S^m (r)$,
			
			\item \label{item:mainTheoremSphere} if $c = 1$, then $\varphi (M)$ is either a small hypersphere $\mathbb S^m (r)$, $r \in (0, 1)$, or the extrinsic product $\mathbb S^{m_1} (r_1) \times \mathbb S^{m_2} (r_2)$, where 
			\begin{equation*}
				r_1^2 \in \left ( 0, \frac 1 {1 + x_2} \right ] \cup \left [ \frac 1 {1 + x_1}, 1 \right ),
			\end{equation*}
			$x_1$ and $x_2$, $x_1 < x_2$, are the solutions of the following second degree polynomial equation
			\begin{equation*}
				m_1 \left ( \frac {2m + 25} {2 (m - 1)} m_1 - 1 \right ) x^2 - \frac {2m + 25} {m - 1} m_1 m_2 x + m_2 \left ( \frac {2m + 25} {2 (m - 1)} m_2 - 1 \right ) = 0.
			\end{equation*}
		\end{enumerate}
	\end{theorem}
	
	\begin{proof}
		Replacing $k = 1$ in \eqref{eq:finalFormulaGenreal}, we obtain
		\begin{equation*}
			\Div \theta_1 \leq |\grad f|^2 \left ( |A|^2 - \frac {m^2 (2m + 25)} {2 (m - 1)} f^2 \right ) - \frac 1 2 f^2 \sum _{i, j = 1} ^m (\lambda_i - \lambda_j)^2 R_{ijij}.
		\end{equation*}
		Integrating and using the hypotheses, we obtain
		\begin{equation*}
			0 \leq \int_M \left ( |\grad f|^2 \left ( |A|^2 - \frac {m^2 (2m + 25)} {2 (m - 1)} f^2 \right ) - \frac 1 2 f^2 \sum _{i, j = 1} ^m (\lambda_i - \lambda_j)^2 R_{ijij} \right ) \ v_g \leq 0.
		\end{equation*}
		Then,
		\begin{equation*}
			|\grad f|^2 \left ( |A|^2 - \frac {m^2 (2m + 25)} {2 (m - 1)} f^2 \right ) = 0
		\end{equation*}
		and 
		\begin{equation*}
			f^2 \sum_{i, j = 1} ^m (\lambda_i - \lambda_j)^2 R_{ijij} = 0,
		\end{equation*}
		on $M$.
		
		In the following, we work on a connected component $V$ of $M_A$ and show that $\grad f = 0$ on $V$. For this, we assume by way of contradiction that $\grad f \neq 0$ at any point of an open subset $U \subset V$. On $U$, there is a point at which $f \neq 0$, otherwise $f = 0$ on $U$ and, thus, $\grad f = 0$ on $U$, contradiction. Thus, eventually restricting $U$, $f \neq 0$ at any point of $U$. Therefore, we have
		\begin{equation*}
			\sum _{i, j = 1} ^m (\lambda_i - \lambda_j)^2 R_{ijij} = 0.
		\end{equation*}
		Since $\Riem^M \geq 0$ on $M$, we deduce that $(\lambda_i - \lambda_j)^2 R_{ijij} = 0$, for any $i, j \in \overline {1, m}$. Using the Gauss equation \eqref{eq:GaussEquation}, we obtain that 
		\begin{equation} \label{eq:relationAtMostTwoDPCProof1}
			(\lambda_i - \lambda_j) (c + \lambda_i \lambda_j) = 0, 
		\end{equation} 
		for any $i, j \in \overline {1, m}$. This implies that on $U$ we have at most two distinct principal curvatures at any point. Since $U$ is a subset of a connected component of $M_A$, we obtain that either we have one distinct principal curvature at any point of $U$, or we have two distinct principal curvatures at any point of $U$. Taking into account that $\grad f \neq 0$ at any point of $U$, we deduce that on $U$ we have exactly two distinct principal curvatures at any point. 
		
		Since $U$ is a subset of a connected component of $M_A$, we can consider a local orthonormal frame field $\{ E_i \}_{i \in \overline {1, m}}$ which diagonalizes the shape operator and we write 
		\begin{equation*}
			\lambda_1 = \ldots = \lambda_{m_1} \quad \text{and} \quad \lambda_{m_1 + 1} = \ldots = \lambda_m.
		\end{equation*}
		Since $\grad f \neq 0$ at any point of $U$ and $M$ is biconservative, we can assume that 
		\begin{equation*}
			E_1 = \frac {\grad f} {|\grad f|}, \quad \lambda_1 = - \frac {mf} 2 \quad \text{and} \quad m_1 = 1.
		\end{equation*}
		Using the fact that $\trace A = mf$, we deduce that 
		\begin{equation*}
			\lambda_2 = \frac {3m} {2 (m - 1)} f.
		\end{equation*}
		Since \eqref{eq:relationAtMostTwoDPCProof1} holds on $U$, we deduce that $c + \lambda_1 \lambda_2 = 0$, which implies that 
		\begin{equation*}
			f^2 = \frac {4 (m - 1) c} {3m^2}. 
		\end{equation*}
		Note that if $c < 0$, then we have a contradiction. If $c \geq 0$, then $f$ is constant on $U$, thus $\grad f = 0$ on $U$, contradiction.
		
		Therefore, $\grad f = 0$ on $V$. This implies that $\grad f = 0$ on $M_A$ and thus, from the density of $M_A$, it holds on $M$. Now, using Theorem \ref{th:NomizuSmythTheorem}, we conclude that $\nabla A = 0$ on $M$.
		
		Now, we check if the hypersurfaces with $\nabla A = 0$ satisfy the hypotheses of the theorem. It is known that when $c \in \{ -1, 0 \}$, the only compact hypersurfaces with $\nabla A = 0$ are the hyperspheres $\mathbb S^m (r)$, $r > 0$. Since they are umbilical, we have $|A|^2 = m f^2$. It is easy to check that $\Riem^M \geq 0$ and \eqref{eq:upperBoundNormAForK=1} is satisfied.
		
		When $c = 1$, the hypersurface $M$ is either a hypersphere $\mathbb S^m (r)$, $r \in (0, 1)$, or an extrinsic product $\mathbb S^{m_1} (r_1) \times \mathbb S^{m_2} (r_2)$, $m_1 + m_2 = m$ and $r_1^2 + r_2^2 = 1$. Again, the hypersphere $\mathbb S^m (r)$ is umbilical and satisfies the hypotheses of the theorem.
		
		Now, we study the extrinsic products $\mathbb S^{m_1} (r_1) \times \mathbb S^{m_2} (r_2)$. The shape operator of this hypersurface is 
		\begin{equation*}
			A (X) = \left ( - \frac {r_2} {r_1} X_1, \frac {r_1} {r_2} X_2 \right ),
		\end{equation*}
		for any $X = (X_1, X_2) \in C\left ( T \mathbb S^{m_1} (r_1) \times \mathbb S^{m_2} (r_2) \right )$, where $X_1 \in C\left ( T \mathbb S^{m_1} (r_1) \right )$ and $X_2 \in C \left ( T \mathbb S^{m_2} (r_2) \right )$. Thus, 
		\begin{equation} \label{eq:normAProducts}
			|A|^2 = \frac {r_2^2} {r_1^2} m_1 + \frac {r_1^2} {r_2^2} m_2
		\end{equation}
		and
		\begin{equation} \label{eq:fProducts}
			f^2 = \frac 1 {m^2} \left ( \frac {r_2^2} {r_1^2} m_1^2 - 2 m_1 m_2 + \frac {r_1^2} {r_2^2} m_2^2 \right ).
		\end{equation}
		Now, we check when 
		\begin{align*}
			& |A|^2 \leq \frac {m^2 (2m + 25)} {2 (m - 1)} f^2 \\
			\Leftrightarrow & \frac {r_2^2} {r_1^2} m_1 + \frac {r_1^2} {r_2^2} m_2 \leq \frac {2m + 25} {2 (m - 1)} \left ( \frac {r_2^2} {r_1^2} m_1^2 - 2 m_1 m_2 + \frac {r_1^2} {r_2^2} m_2^2 \right ) \\
			\Leftrightarrow & m_1 \left ( \frac {2m + 25} {2 (m - 1)} m_1 - 1 \right ) \left ( \frac {r_2^2} {r_1^2} \right )^2 - \frac {2m + 25} {m - 1} m_1 m_2 \frac {r_2^2} {r_1^2} + m_2 \left ( \frac {2m + 25} {2 (m - 1)} m_2 - 1 \right ) \geq 0.
		\end{align*}
		Studying the coefficients and the discriminant of the second degree polynomial attached to the above inequality, we deduce that its roots are positive real numbers. Since, the leading coefficient is positive in our case, we deduce that 
		\begin{equation*}
			\frac {r_2^2} {r_1^2} \in (0, x_1] \cup [x_2, \infty),
		\end{equation*}
		where $x_1$ and $x_2$ are the roots mentioned earlier.
		
		Using the fact that $r_2^2 = 1 - r_1^2$, we obtain that 
		\begin{equation*}
			r_1^2 \in \left ( 0, \frac 1 {1 + x_2} \right ] \cup \left [\frac 1 {1 + x_1}, 1 \right ).
		\end{equation*}
	\end{proof}
	
	\begin{remark}
		We note that the minimal extrinsic products are not included in the conclusion of Theorem \ref{th:mainTheorem}. Indeed, the extrinsic product $\mathbb S^{m_1} (r_1) \times \mathbb S^{m_2} (r_2)$ is minimal if and only if $r_1^2 = m_1 / m$, i.e. $r_2^2 / r_1^2 = m_2 / m_1$. Replacing this in \eqref{eq:normAProducts} and \eqref{eq:fProducts}, we obtain $|A|^2 = m$ and $f = 0$, so \eqref{eq:upperBoundNormAForK=1} does not hold. Therefore, $m_2 / m_1 \in (x_1, x_2)$ and $r_1^2 = m_1 / m$ does not belong to the intervals from Theorem \ref{th:mainTheorem}. This explains why, unlike in Theorem \ref{th:FetcuOniciuc2022}, the condition $r_1 \neq \sqrt {m_1 / m}$ does not appear in our result. Moreover, the upper bound \eqref{eq:upperBoundNormAForK=1} is attained. Indeed, if $r_1^2 = 1 / (1 + x_1)$ or $r_1^2 = 1 / (1 + x_2)$, then $r_2^2 / r_1^2 = x_1$ or $r_2^2 / r_1^2 = x_2$, respectively, and, replacing in \eqref{eq:normAProducts} and \eqref{eq:fProducts}, we obtain equality in \eqref{eq:upperBoundNormAForK=1} at any point of the extrinsic product. From the above discussion, these hypersurfaces are not minimal. However, since they are CMC, they do not show that the upper bound \eqref{eq:upperBoundNormAForK=1} cannot be improved. The same holds for the upper bound \eqref{eq:upperBoundNormAForK0} and the extrinsic products from Theorem \ref{th:mainTheoremFneq0}.
	\end{remark}
	
	In the proof of Theorem \ref{th:mainTheorem}, the term which contains the upper bound of $|A|^2$ is multiplied by $|\grad f|^2$. Therefore, it is enough to assume that \eqref{eq:upperBoundNormAForK=1} is satisfied only at the points of $M$ where $\grad f \neq 0$, and we obtain the following result. The condition on $r_1$ from Theorem \ref{th:mainTheorem} is no longer obtained, since the extrinsic products are CMC.
	
	\begin{corollary} \label{cor:upperBoundOnlyWhereGradFNeq0}
		Let $\varphi : M^m \to N^{m+1} (c)$ be a compact non-minimal biconservative hypersurface with $m > 10$. If $\Riem^M \geq 0$ and \eqref{eq:upperBoundNormAForK=1} holds at the points of $M$ where $\grad f \neq 0$, then $\nabla A = 0$ and $\varphi (M)$ is one of the hypersurfaces in Theorem \ref{th:mainTheorem}, where the condition on $r_1$ in item \ref{item:mainTheoremSphere} is replaced by $r_1 \neq \sqrt {m_1 / m}$.
	\end{corollary}
	
	\begin{corollary} \label{cor:atMostTwoDPC}
		Let $\varphi : M^m \to N^{m+1} (c)$ be a compact non-minimal biconservative hypersurface with $m > 10$ and $\Riem^M \geq 0$. If $M$ has at most two distinct principal curvatures at the points where $\grad f \neq 0$, then $\nabla A = 0$ and $\varphi (M)$ is one of the hypersurfaces in Theorem \ref{th:mainTheorem}, where the condition on $r_1$ in item \ref{item:mainTheoremSphere} is replaced by $r_1 \neq \sqrt {m_1 / m}$.
	\end{corollary}
	
	\begin{proof}
		As in Corollary \ref{cor:upperBoundOnlyWhereGradFNeq0}, it is enough to show that
		\begin{equation} \label{eq:inequalityGradFNormA}
			|\grad f|^2 \left ( |A|^2 - \frac {m^2 (2m + 25)} {2 (m - 1)} f^2 \right ) \leq 0
		\end{equation}
		on $M$.
		
		To prove this, let $V$ be a connected component of the open and dense set $M_A$ of all points at which the number of distinct principal curvatures is locally constant. Consider an arbitrary point $p \in V$.
		
		If $\grad f = 0$ at $p$, then \eqref{eq:inequalityGradFNormA} holds at $p$.
		
		If $\grad f \neq 0$ at $p$, then, since $V$ is open, there is an open neighbourhood $U \subset V$ of $p$ such that $\grad f \neq 0$ at any point of $U$. Since $U \subset V$, from the hypotheses, we deduce that on $U$ we have two distinct principal curvatures. Thus, using Lemma \ref{th:LemmaBiconservativeHypersurfacesConnectionForms}, we obtain
		\begin{equation*}
			\lambda_1 = - \frac {m} 2 f \quad \text{and} \quad \lambda_2 = \ldots = \lambda_m = \frac {3m} {2 (m - 1)} f, \quad \text{on } U.
		\end{equation*}
		Using this, we compute
		\begin{equation*}
			|A|^2 = \frac {m^2 (m + 8)} {4 (m - 1)} f^2 \leq \frac {m^2 (2m + 25)} {2 (m - 1)} f^2.
		\end{equation*}
		Thus, \eqref{eq:inequalityGradFNormA} holds at $p$.
		
		Therefore, \eqref{eq:inequalityGradFNormA} holds on $V$, and thus on $M_A$. Using the density of $M_A$ and the continuity of the left-hand side of \eqref{eq:inequalityGradFNormA}, we obtain that \eqref{eq:inequalityGradFNormA} holds on $M$. Following the proof of Theorem \ref{th:mainTheorem}, we obtain the conclusion.
	\end{proof}
	
	We have seen that, when $m > 10$, the right-hand side of \eqref{eq:finalFormulaGenreal} is non-positive if 
	\begin{equation} \label{eq:upperBoundNormAForK0}
		|A|^2 \leq \frac {m^2} {4 (m - 1)} b_{k_0} f^2 = \frac {m^2 (2m^2 + 5m - 34)} {2 (m - 1) (m - 10)} f^2.
	\end{equation}
	This upper bound cannot be reached by any positive integer $k$. However, inequality \eqref{eq:finalFormulaGenreal} can be extended to any real number $k$, if we impose that $f \neq 0$ at any point of $M$. Indeed, eventually changing the orientation of $M$, we can suppose that $f > 0$ on $M$ and the proof of \eqref{eq:finalFormulaGenreal} holds for any real number $k$. 
	
	Following the same steps as in the proof of Theorem \ref{th:mainTheorem}, we obtain
	\begin{theorem} \label{th:mainTheoremFneq0}
		Let $\varphi : M^m \to N^{m+1} (c)$ be a compact biconservative hypersurface with $m > 10$. If $\Riem^M \geq 0$, $f \neq 0$ at any point of $M$ and \eqref{eq:upperBoundNormAForK0} holds on $M$, then $\nabla A = 0$. Moreover,
		\begin{enumerate} [label = \alph*)]
			\item if $c = -1$, then $\varphi (M)$ is a hypersphere $\mathbb S^m (r)$ of radius $r > 0$,
			
			\item if $c = 0$, then $\varphi (M)$ is a hypersphere $\mathbb S^m (r)$,
			
			\item \label{item:mainTheoremFneq0Sphere} if $c = 1$, then $\varphi (M)$ is either a small hypersphere $\mathbb S^m (r)$, $r \in (0, 1)$, or the extrinsic product $\mathbb S^{m_1} (r_1) \times \mathbb S^{m_2} (r_2)$, where 
			\begin{equation*}
				r_1^2 \in \left ( 0, \frac 1 {1 + x_2} \right ] \cup \left [ \frac 1 {1 + x_1}, 1 \right ),
			\end{equation*}
			$x_1$ and $x_2$, $x_1 < x_2$, are the solutions of the following second degree polynomial equation
			\begin{align*}
				& m_1 \left ( \frac {2m^2 + 5m -34} {2 (m - 1) (m - 10)} m_1 - 1 \right ) x^2 - \frac {2m^2 + 5m - 34} {(m - 1) (m - 10)} m_1 m_2 x \\
				& + m_2 \left ( \frac {2m^2 + 5m - 34} {2 (m - 1) (m - 10)} m_2 - 1 \right ) = 0.
			\end{align*}
		\end{enumerate}
	\end{theorem}
	
	As in the case of Theorem \ref{th:mainTheorem}, in the proof of Theorem \ref{th:mainTheoremFneq0} the term which contains the upper bound of $|A|^2$ is multiplied by $|\grad f|^2$. Therefore, we obtain the following result.
	
	\begin{corollary} \label{cor:upperBoundK0OnlyWhereGradFNeq0}
		Let $\varphi : M^m \to N^{m+1} (c)$ be a compact biconservative hypersurface with $m > 10$. If $\Riem^M \geq 0$, $f \neq 0$ at any point of $M$ and \eqref{eq:upperBoundNormAForK0} holds at the points of $M$ where $\grad f \neq 0$, then $\nabla A = 0$ and $\varphi (M)$ is one of the hypersurfaces in Theorem \ref{th:mainTheoremFneq0}, where the condition on $r_1$ in item \ref{item:mainTheoremFneq0Sphere} is replaced by $r_1 \neq \sqrt {m_1 / m}$.
	\end{corollary}
	
	We note that if the left-hand side of \eqref{eq:finalFormulaGenreal} would vanish, then we could follow the proofs of Theorems \ref{th:mainTheorem} and \ref{th:mainTheoremFneq0} without the hypothesis of compactness. In order to do this, we compute $\theta_1$. Using \eqref{eq:formulaTheta} and the fact that $\varphi$ is biconservative, we obtain
	\begin{align*}
		\theta_1 =& - \left ( \Div T_1^2 \right ) ^\sharp \\
				 =& - \left ( \Div \left ( f^2 A^2 - \frac m 2 f^3 A + \frac {m^2} {16} f^4 \id \right ) \right )^\sharp \\
				 =& - \left ( A^2 \left (\grad f^2 \right ) + f^2 \left ( \Div A^2 \right )^\sharp - \frac m 2 A \left ( \grad f^3 \right ) - \frac m 2 f^3 \left ( \Div A \right )^\sharp + \frac {m^2} {16} \grad f^4 \right ) \\
				 =& - \left ( \frac {m^2} 2 f^3 \grad f + f^2 \left ( \Div A^2 \right )^\sharp + \frac {3m^2} 4 f^3 \grad f - \frac {m^2} 2 f^3 \grad f + \frac {m^2} 4 f^3 \grad f \right ) \\
				 =& - \left ( m^2 f^3 \grad f + f^2 \left ( \Div A^2 \right )^\sharp \right ).
	\end{align*}
	Therefore, using item \ref{item:divA2} of Lemma \ref{th:PropertiesNablaA} and \eqref{eq:BiconservativityCharacterization}, we obtain
	\begin{equation*}
		\theta_1 = - \frac {f^2} 2 \grad \left ( |A|^2 + \frac {m^2} 2 f^2 \right ) .
	\end{equation*}
	Similarly, 
	\begin{equation*}
		\theta_{k_0} = - \frac {f^{2k_0}} 2 \grad \left ( |A|^2 + \frac {m^2} 2 f^2 \right ).
	\end{equation*}
	We note that if $|A|^2 + m^2 f^2 / 2$ is constant, then $\Div \theta_1 = 0$ and $\Div \theta_{k_0} = 0$, that is the left-hand side of \eqref{eq:mainTensorialFormula} vanishes in both cases. Using this remark and the proofs of Theorems \ref{th:mainTheorem} and \ref{th:mainTheoremFneq0}, we obtain the following two results. 
	
	\begin{theorem} \label{th:mainTheoremWithoutCompactness}
		Let $\varphi : M^m \to N^{m+1} (c)$ be a non-minimal biconservative hypersurface with $m > 10$. If $\Riem^M \geq 0$, $|A|^2 + m^2 f^2 / 2$ is constant and \eqref{eq:upperBoundNormAForK=1} holds on $M$, then $\nabla A = 0$. Moreover,
		\begin{enumerate} [label = \alph*)]
			\item if $c = -1$, then $\varphi (M)$ is either an open part of a hypersphere $\mathbb S^m (r)$ of radius $r > 0$, or an open part of a horosphere, or an open part of the standard product $\mathbb S^{m - 1} (r) \times \mathbb H^1 \left ( \sqrt {1 + r^2} \right )$, $r > 0$, i.e. a tube around a geodesic of $\mathbb H^{m+1}$,

			\item if $c = 0$, then $\varphi (M)$ is either an open part of a hypersphere $\mathbb S^m (r)$, or an open part of the generalized cylinder $\mathbb S^{m_1} (r) \times \mathbb R^{m_2}$, where $r > 0$, $m_1 + m_2 = m$, $m_1 \geq 1$ and $m_2 \geq 1$,
			
			\item if $c = 1$, then $\varphi (M)$ is either an open part of a small hypersphere $\mathbb S^m (r)$, $r \in (0, 1)$, or an open part of the extrinsic product $\mathbb S^{m_1} (r_1) \times \mathbb S^{m_2} (r_2)$, where 
			\begin{equation*}
				r_1^2 \in \left ( 0, \frac 1 {1 + x_2} \right ] \cup \left [ \frac 1 {1 + x_1}, 1 \right ),
			\end{equation*}
			$x_1$ and $x_2$, $x_1 < x_2$, are the solutions of the following second degree polynomial equation
			\begin{equation*}
				m_1 \left ( \frac {2m + 25} {2 (m - 1)} m_1 - 1 \right ) x^2 - \frac {2m + 25} {m - 1} m_1 m_2 x + m_2 \left ( \frac {2m + 25} {2 (m - 1)} m_2 - 1 \right ) = 0.
			\end{equation*}
		\end{enumerate}
	\end{theorem}
	
	\begin{proof}
		Since $|A|^2 + m^2 f^2 / 2$ is constant, we have $\theta_1 = 0$ and thus $\Div \theta_1 = 0$. Therefore, the left-hand side of \eqref{eq:finalFormulaGenreal} vanishes for $k = 1$ and, following the proof of Theorem \ref{th:mainTheorem}, we obtain that $f$ is constant. Since $|A|^2 + m^2 f^2 / 2$ and $f$ are constant, it follows that $|A|^2$ is constant too. Now, the conclusion follows from Remark \ref{rem:localNomizuSmyth}, taking into account that $M$ is non-minimal. Using Theorem 4 of \cite{Lawson1969}, we obtain the classification of the hypersurfaces.
	\end{proof}
	
	\begin{theorem} \label{th:mainTheoremFneq0WithoutCompactness}
		Let $\varphi : M^m \to N^{m+1} (c)$ be a biconservative hypersurface with $m > 10$. If $\Riem^M \geq 0$, $f \neq 0$ at any point of $M$, $|A|^2 + m^2 f^2 / 2$ is constant and \eqref{eq:upperBoundNormAForK0} holds on $M$, then $\nabla A = 0$. Moreover,
		\begin{enumerate} [label = \alph*)]
			\item if $c = -1$, then $\varphi (M)$ is either an open part of a hypersphere $\mathbb S^m (r)$ of radius $r > 0$, or an open part of a horosphere, or an open part of the standard product $\mathbb S^{m - 1} (r) \times \mathbb H^1 \left ( \sqrt {1 + r^2} \right )$, $r > 0$, i.e. a tube around a geodesic of $\mathbb H^{m+1}$,

			\item if $c = 0$, then $\varphi (M)$ is either an open part of a hypersphere $\mathbb S^m (r)$, or an open part of the generalized cylinder $\mathbb S^{m_1} (r) \times \mathbb R^{m_2}$, where $r > 0$, $m_1 + m_2 = m$, $m_1 \geq 1$ and $m_2 \geq 1$,
			
			\item if $c = 1$, then $\varphi (M)$ is either an open part of a small hypersphere $\mathbb S^m (r)$, $r \in (0, 1)$, or an open part of the extrinsic product $\mathbb S^{m_1} (r_1) \times \mathbb S^{m_2} (r_2)$, where 
			\begin{equation*}
				r_1^2 \in \left ( 0, \frac 1 {1 + x_2} \right ] \cup \left [ \frac 1 {1 + x_1}, 1 \right ),
			\end{equation*}
			$x_1$ and $x_2$, $x_1 < x_2$, are the solutions of the following second degree polynomial equation
			\begin{align*}
				& m_1 \left ( \frac {2m^2 + 5m -34} {2 (m - 1) (m - 10)} m_1 - 1 \right ) x^2 - \frac {2m^2 + 5m - 34} {(m - 1) (m - 10)} m_1 m_2 x \\
				& + m_2 \left ( \frac {2m^2 + 5m - 34} {2 (m - 1) (m - 10)} m_2 - 1 \right ) = 0.
			\end{align*}
		\end{enumerate}
	\end{theorem}
	
	\vspace{15pt}
	
	\noindent \textbf{Acknowledgments.} We are grateful to Cezar Oniciuc for his valuable suggestions and discussions.
	
	\bibliographystyle{abbrv}
	\bibliography{Bibliography.bib}
\end{document}